\documentclass{LMCS}

\def\dOi{10(3:1)2014}
\lmcsheading%
{\dOi}
{1--23}
{}
{}
{Mar.~30, 2013}
{Aug.~\phantom07, 2014}
{}

\ACMCCS{[{\bf Theory of computation}]: Logic}

\usepackage{hyperref}
\usepackage{proof}

\newcommand{\limplies}{\mathbin{\rightarrow}}

\newcommand{\metaimplies}{\mathrel{\Rightarrow}}

\newcommand{\Pos}{{\mathcal{V}^+}}
\newcommand{\Nsp}{{\mathcal{V}^+_{ns}}}
\newcommand{\Neg}{{\mathcal{V}^-}}

\begin{document}

\title[decidability of variables]{Classical propositional logic and decidability of variables in
  intuitionistic propositional logic}

\author[H.~Ishihara]{Hajime Ishihara}
\address{
School of Information Science, Japan Advanced Institute of Science and Technology, Nomi, Ishikawa 923-1292, Japan
}	%required
\email{
ishihara@jaist.ac.jp
}  %optional
%\thanks{thanks 1, optional.}	%optional

%% required for running head on odd and even pages, use suitable
%% abbreviations in case of long titles and many authors:

%% mandatory lists of keywords and classifications:
%%%%%%%%%%%%%%%%%%%%%%%%%%%%%%%%%%%%%%%%%%%%%%%%%%%%%%%%%%%%%%%%%%%%%%%%%%%

%% the abstract has to PRECEED the command \maketitle:
%% be sure not to issue the \maketitle command twice!

\begin{abstract}
We improve the answer to the question: {\it what set of
  excluded middles for propositional variables in a formula suffices to prove the formula in intuitionistic propositional logic whenever it is provable in classical propositional logic}.
\end{abstract}

\maketitle

\section{Introduction}

Let $\vdash_c$ and $\vdash_i$ denote derivability in classical and
intuitionistic propositional logic, respectively.
Then it is known that if $\vdash_c A$, then $\Pi_{\mathcal{V}(A)}
\vdash_i A$, where $\mathcal{V}(A)$ is the set of propositional
variables in a formula $A$ and $\Pi_V = \{ p \lor \lnot p
\mid p \in V \}$ for a set $V$ of propositional variables;
see, for example, \cite[appendix]{Fu97}, and \cite[p. 27]{NvP01} which was
originally given in \cite{vP98}.

In this note, we consider a problem:
{\it what set $V$ of propositional variables suffices for $\Pi_V,
\Gamma \vdash_i A$ whenever $\Gamma \vdash_c A$}, and show, employing
a technique in \cite{Is00,Is12}, that
$ V = (\Neg(\Gamma) \cup \Pos(A)) \cap (\Nsp(\Gamma) \cup \Neg(A)) $
suffices, where $\Pos$, $\Neg$ and $\Nsp$ are the sets of
propositional variables occurring positively, negatively and
non-strictly positively, respectively (precise definitions will be
given in the next section).
For example, since 
$ (p \limplies q) \limplies p \vdash_c p , $
we have
\[ p \lor \lnot p, (p \limplies q) \limplies p \vdash_i p \]
and, since
$ p \limplies q \lor r \vdash_c (p \limplies q) \lor (p \limplies r) , $
we have 
\[ p \lor \lnot p, p \limplies q \lor r \vdash_i (p \limplies q) \lor
(p \limplies r) . \]
\newpage
\section{Preliminaries}

We refer to Troelstra and Schwichtenberg \cite{TS96} for the
necessary background on sequent calculi; see also Negri and von
Plato \cite{NvP01}.
We use the standard language of propositional logic containing
$\land$, $\lor$, $\limplies$ and $\bot$ as primitive logical
operators, and introduce the abbreviation $\lnot A \equiv A \limplies
\bot$.
We define {\it positive}, {\it strictly positive} and {\it negative
occurrence} of a formula in the usual way (see \cite[1.1.3]{TS96}
or \cite[3.9,3.11,3.23]{TD88} for details).
The sets $\Pos(A)$ and $\Neg(A)$ of propositional variables occurring
positively and negatively, respectively, in a formula $A$ are
simultaneously defined by
\begin{eqnarray*}
\Pos(p) & = & \{p\} , \quad \Pos(\bot) = \emptyset, \\
\Pos(A \land B) & = & \Pos(A \lor B) = \Pos(A) \cup \Pos(B) , \\ 
\Pos(A \limplies B) & = & \Neg(A) \cup \Pos(B) , \\
\Neg(p) & = & \Neg(\bot) = \emptyset , \\
\Neg(A \land B) & = & \Neg(A \lor B) = \Neg(A) \cup \Neg(B) , \\ 
\Neg(A \limplies B) & = & \Pos(A) \cup \Neg(B) .
\end{eqnarray*}
The set $\Nsp(A)$ of propositional variables occurring
non-strictly positively in a formula $A$ is defined by
\begin{eqnarray*}
\Nsp(p) & = & \Nsp(\bot) = \emptyset , \\
\Nsp(A \land B) & = & \Nsp(A \lor B) = \Nsp(A) \cup \Nsp(B) , \\
\Nsp(A \limplies B) & = & \Neg(A) \cup \Nsp(B) .
\end{eqnarray*}
We extend $\Pos$ to a finite multiset $\Gamma$ of
formulas by $\Pos(\Gamma) = \bigcup_{A \in \Gamma} \Pos(A)$.
$\Neg(\Gamma)$ and $\Nsp(\Gamma)$ are defined similarly.

The sequent calculus ${\bf G3cp}$ is specified by the following
axioms and rules:
$$
\begin{array}{cc}
 p, \Gamma \metaimplies \Delta, p \quad {\rm Ax}
 &
 \bot, \Gamma \metaimplies \Delta \quad {\rm L}\bot \\[2mm]
 \infer[{\rm L}\land]{A \land B, \Gamma \metaimplies \Delta }{
  A, B, \Gamma \metaimplies \Delta
 }
 &
 \infer[{\rm R}\land]{\Gamma \metaimplies \Delta, A \land B}{
  \Gamma \metaimplies \Delta, A
  &
  \Gamma \metaimplies \Delta, B
 } \\[2mm]
 \infer[{\rm L}\lor]{A \lor B, \Gamma \metaimplies \Delta}{
  A, \Gamma \metaimplies \Delta
  &
  B, \Gamma \metaimplies \Delta
 }
 &
 \infer[{\rm R}\lor]{\Gamma \metaimplies \Delta, A \lor B}{
  \Gamma \metaimplies \Delta, A, B
 } \\[2mm]
 \infer[{\rm L}\limplies]{A \limplies B, \Gamma \metaimplies \Delta}{
  \Gamma \metaimplies \Delta, A
  &
  B, \Gamma \metaimplies \Delta
 }
 &
 \infer[{\rm R}\limplies]{\Gamma \metaimplies \Delta, A \limplies B}{
  A, \Gamma \metaimplies \Delta, B
 }
\end{array}
$$
where in ${\rm Ax}$, $p$ is a propositional variable.

The intuitionistic version ${\bf G3ip}$ of ${\bf G3cp}$ has the following form:
$$
\begin{array}{cc}
 p, \Gamma \metaimplies p \quad {\rm Ax}
 &
 \bot, \Gamma \metaimplies A \quad {\rm L}\bot \\[2mm]
 \infer[{\rm L}\land]{A \land B, \Gamma \metaimplies C}{
  A, B, \Gamma \metaimplies C
 }
 &
 \infer[{\rm R}\land]{\Gamma \metaimplies A \land B}{
  \Gamma \metaimplies A
  &
  \Gamma \metaimplies B
 } \\[2mm]
 \infer[{\rm L}\lor]{A \lor B, \Gamma \metaimplies C}{
  A, \Gamma \metaimplies C
  &
  B, \Gamma \metaimplies C
 }
 &
 \infer[{\rm R}\lor_1]{\Gamma \metaimplies A \lor B}{
  \Gamma \metaimplies A
 } \quad
 \infer[{\rm R}\lor_2]{\Gamma \metaimplies A \lor B}{
  \Gamma \metaimplies B
 } \\[2mm]
 \infer[{\rm L}\limplies]{A \limplies B, \Gamma \metaimplies C}{
  A \limplies B, \Gamma \metaimplies A
  &
  B, \Gamma \metaimplies C
 }
 &
 \infer[{\rm R}\limplies]{\Gamma \metaimplies A \limplies B}{
  A, \Gamma \metaimplies B
 }
\end{array}
$$
where in ${\rm Ax}$, $p$ is a propositional variable.

Note that having the present sequent calculus formulation (${\rm Ax}$
with a propositional variable $p$ instead of a formula $A$) allows for
an easy treatment of the Basis case in the proof of the main result
below.

The structural rules (weakening, contraction and cut) are admissible 
in ${\bf G3cp}$ and in ${\bf G3ip}$; see \cite[3.4.3,3.4.5,4.1.2]{TS96}.
Those structural rules are formulated in ${\bf G3ip}$ as follows:
$$
\begin{array}{cc}
\infer[{\rm LW}]{
  \Gamma, \Delta \metaimplies C
}{
  \Gamma \metaimplies C
}
&
\infer[{\rm LC}]{
  A, \Gamma \metaimplies C
}{
  A, A, \Gamma \metaimplies C
}
\end{array}
$$
$$
\infer[{\rm Cut}]{
  \Gamma, \Gamma' \metaimplies C
}{
  \Gamma \metaimplies A
  &
  A, \Gamma' \metaimplies C
} .
$$

We write $\vdash_c \Gamma \metaimplies \Delta$ and $\vdash_i \Gamma
\metaimplies A$ for derivability of sequents $\Gamma \metaimplies
\Delta$ and $\Gamma \metaimplies A$ in ${\bf G3cp}$ and in ${\bf
G3ip}$, respectively.

We introduce the symbol ``$\ast$'' as a special proposition letter
(a {\it place holder}) and an abbreviation $\lnot_\ast A \equiv A
\limplies \ast$.
It is straightforward to see that if $\vdash_i \Gamma \metaimplies
A$ then $\vdash_i \Gamma, \lnot_\ast A \metaimplies \ast$; if
$\vdash_i \Gamma, \lnot_\ast \lnot_\ast A \metaimplies \ast$ then
$\vdash_i \Gamma \metaimplies \lnot_\ast A$, and $\vdash_i
\Gamma, A \metaimplies \ast$ if and only if $\vdash_i \Gamma
\metaimplies \lnot_\ast A$.
From the latter and the former results, it is trivial to conclude that
if $\vdash_i \Gamma, A \metaimplies \ast$ then $\vdash_i \Gamma,
\lnot_\ast \lnot_\ast A \metaimplies \ast$, and $\vdash_i \Gamma,
\lnot_\ast A \metaimplies \ast$ if and only if $\vdash_i \Gamma
\metaimplies \lnot_\ast \lnot_\ast A$.

We have the following lemma for the logical operators and the
operators $\lnot$ and $\lnot_\ast$.

\begin{lem}\label{basics}\hfill
\begin{enumerate}
\item $\vdash_i \Gamma, p \lor \lnot p, \lnot_\ast \lnot p, \lnot_\ast
  p \metaimplies \ast$,\label{EM}
\item $\vdash_i \Gamma, \lnot_\ast \lnot \bot \metaimplies \ast$,\label{Lbot}
\item $\vdash_i \lnot_\ast \lnot (D \land D') \metaimplies
  \lnot_\ast \lnot D  \land \lnot_\ast \lnot D'$,\label{Lland}
\item $\vdash_i \lnot_\ast \lnot_\ast S \land \lnot_\ast \lnot_\ast S'
  \metaimplies \lnot_\ast   \lnot_\ast (S \land S')$,\label{Rland}
\item $\vdash_i \lnot_\ast \lnot (D \lor D') \metaimplies
  \lnot_\ast \lnot_\ast (\lnot_\ast \lnot D \lor \lnot_\ast \lnot
  D')$,\label{Llor}
\item $\vdash_i \lnot_\ast (\lnot_\ast S \land \lnot_\ast S')
  \metaimplies \lnot_\ast \lnot_\ast (S \lor S')$,\label{Rlor}
\item $\vdash_i \lnot_\ast \lnot (S \limplies B) \metaimplies
  \lnot_\ast \lnot_\ast S \limplies \lnot_\ast \lnot
  B$,\label{Llimplies2}
\item $\vdash_i S \limplies B \metaimplies \lnot_\ast
  \lnot_\ast S \limplies \lnot_\ast \lnot_\ast B$,\label{Llimplies1}
\item $\vdash_i \lnot_\ast \lnot A \limplies \lnot_\ast
  \lnot_\ast S \metaimplies \lnot_\ast \lnot_\ast (A \limplies
  S)$.\label{Rlimplies}
\end{enumerate}
\end{lem}
\proof
Easy exercise.
\qed

Let $A[\ast/C]$ denote the result of substituting a formula $C$ for
each occurrence of $\ast$ in a formula $A$, and, for a finite multiset
$\Gamma \equiv A_1, \ldots, A_n$, let $\Gamma[\ast/C]$ denote the
multiset $A_1[\ast/C], \ldots, A_n[\ast/C]$.

\begin{lem}\label{subst}
If
$ \vdash_i \Gamma \metaimplies A , $
then
$ \vdash_i \Gamma[\ast/C] \metaimplies A[\ast/C] . $
\end{lem}
\proof
By induction on the depth of a deduction
$ \vdash_i \Gamma \metaimplies A . $
\qed

\section{The main result}

If ``$c$'' is an operator, such as $\lnot$ and $\lnot_\ast$, and
$\Gamma \equiv A_1, \ldots, A_n$ is a finite multiset of formulas,
then we write $c \Gamma$ for the multiset $c A_1, \ldots, c A_n$.

\begin{prop}\label{main}
If 
$ \vdash_c \Gamma, \Delta \metaimplies \Sigma , $
then
$ \vdash_i \Pi_V, \Gamma, \lnot_\ast \lnot \Delta, \lnot_\ast \Sigma
\metaimplies \ast , $
where $V$ is a set of propositional variables containing
$ (\Neg(\Gamma, \Delta) \cup \Pos(\Sigma)) \cap (\Nsp(\Gamma)
\cup \Pos(\Delta) \cup \Neg(\Sigma)) . $
\end{prop}
\proof
Let $V$ be a set of propositional variables containing
$ (\Neg(\Gamma, \Delta) \cup \Pos(\Sigma)) \cap (\Nsp(\Gamma)
\cup \Pos(\Delta) \cup \Neg(\Sigma)) , $
and we proceed by induction on the depth of a deduction of 
$ \vdash_c \Gamma, \Delta \metaimplies \Sigma . $

\noindent
{\it Basis.}
If the deduction is an instance of ${\rm Ax}$, then it must be either
of the form
$ p, \Gamma', \Delta \metaimplies \Sigma', p , $
or of the form
$ \Gamma, p, \Delta' \metaimplies \Sigma', p . $
In the former case, we have
\[ \vdash_i \Pi_V, p, \Gamma', \lnot_\ast \lnot \Delta, \lnot_\ast
\Sigma',  \lnot_\ast p \metaimplies \ast \]
and, in the latter case, since
\[ p \in (\Neg(\Gamma, p, \Delta') \cup \Pos(\Sigma', p)) \cap
(\Nsp(\Gamma) \cup \Pos(p, \Delta') \cup \Neg(\Sigma', p)) \subseteq V
, \]
we have
\[ \vdash_i \Pi_V, \Gamma, \lnot_\ast \lnot p, \lnot_\ast \lnot \Delta',
\lnot_\ast \Sigma', \lnot_\ast p \metaimplies \ast \]
by Lemma \ref{basics} (\ref{EM}).
If the deduction is an instance of ${\rm L}\bot$, then it must be
either of the form
$ \bot, \Gamma', \Delta \metaimplies \Sigma , $
or of the form
$ \Gamma, \bot, \Delta' \metaimplies \Sigma . $
In the former case, we have
\[ \vdash_i \Pi_V, \bot, \Gamma', \lnot_\ast \lnot \Delta, \lnot_\ast
\Sigma \metaimplies \ast \]
and, in the latter case, we have
\[ \vdash_i \Pi_V, \Gamma, \lnot_\ast \lnot \bot, \lnot_\ast \lnot
\Delta', \lnot_\ast \Sigma \metaimplies \ast \]
by Lemma \ref{basics} (\ref{Lbot}).

\noindent
{\it Induction step.}
For the induction step, we distinguish the cases:
(A) the last rule applied is an L-rule and the principal formula is in
$\Delta$,
(B) the last rule applied is an L-rule and the principal formula is in
$\Gamma$, and
(C) the last rule applied is an R-rule.

\noindent
{\it Case A.}
The last rule applied is an ${\rm L}$-rule, and the principal
formula is in $\Delta$.

\noindent
{\it Case A1.}
The last rule applied is ${\rm L}\land$.
Then the derivation ends with
$$
 \infer[{\rm L}\land]{
  \Gamma, D \land D', \Delta' \metaimplies \Sigma
  }{
  \Gamma, D, D', \Delta' \metaimplies \Sigma
 } .
$$
Since
\begin{eqnarray*}
\lefteqn{(\Neg(\Gamma, D, D', \Delta') \cup \Pos(\Sigma)) \cap
  (\Nsp(\Gamma) \cup \Pos(D, D', \Delta') \cup \Neg(\Sigma)) =} \\
& &
(\Neg(\Gamma, D \land D', \Delta') \cup \Pos(\Sigma)) \cap (\Nsp(\Gamma)
\cup \Pos(D \land D', \Delta') \cup \Neg(\Sigma)) \subseteq V ,
\end{eqnarray*}
we have
\[ \vdash_i \Pi_V, \Gamma, \lnot_\ast \lnot D, \lnot_\ast \lnot D',
\lnot_\ast \lnot \Delta', \lnot_\ast \Sigma \metaimplies \ast \]
by the induction hypothesis, and hence
\[ \vdash_i \Pi_V, \Gamma, \lnot_\ast \lnot D \land \lnot_\ast \lnot D',
\lnot_\ast \lnot \Delta', \lnot_\ast \Sigma \metaimplies \ast \]
by ${\rm L}\land$.
Therefore
$ \vdash_i \Pi_V, \Gamma, \lnot_\ast \lnot (D \land D'), \lnot_\ast \lnot
\Delta', \lnot_\ast \Sigma \metaimplies \ast , $
by ${\rm Cut}$ with Lemma \ref{basics} (\ref{Lland}).

\noindent
{\it Case A2.}
The last rule applied is ${\rm L}\lor$.
Then the derivation ends with
$$
 \infer[{\rm L}\lor]{
  \Gamma, D \lor D', \Delta' \metaimplies \Sigma
  }{
  \Gamma, D, \Delta' \metaimplies \Sigma
  &
  \Gamma, D', \Delta' \metaimplies \Sigma
 } .
$$
Since
$ (\Neg(\Gamma, D, \Delta') \cup \Pos(\Sigma)) \cap (\Nsp(\Gamma) \cup
\Pos(D, \Delta') \cup \Neg(\Sigma)) \subseteq V $
and
$ (\Neg(\Gamma, D', \Delta') \cup \Pos(\Sigma)) \cap (\Nsp(\Gamma)
\cup \Pos(D', \Delta') \cup \Neg(\Sigma)) \subseteq V , $
we have
\[ \vdash_i \Pi_V, \Gamma, \lnot_\ast \lnot D, \lnot_\ast \lnot \Delta',
\lnot_\ast \Sigma \metaimplies \ast
\quad\mbox{and}\quad
\vdash_i \Pi_V, \Gamma, \lnot_\ast \lnot D', \lnot_\ast \lnot \Delta',
\lnot_\ast \Sigma \metaimplies \ast \]
by the induction hypothesis, and hence
\[ \vdash_i \Pi_V, \Gamma, \lnot_\ast \lnot D \lor \lnot_\ast \lnot D',
\lnot_\ast \lnot \Delta', \lnot_\ast \Sigma \metaimplies \ast \]
by ${\rm L}\lor$.
Therefore
\[ \vdash_i \Pi_V, \Gamma, \lnot_\ast \lnot_\ast (\lnot_\ast \lnot D \lor
\lnot_\ast \lnot D'), \lnot_\ast \lnot \Delta', \lnot_\ast \Sigma
\metaimplies \ast \]
and so
$ \vdash_i \Pi_V, \Gamma, \lnot_\ast \lnot (D \lor D'), \lnot_\ast \lnot
\Delta', \lnot_\ast \Sigma \metaimplies \ast , $
by ${\rm Cut}$ with Lemma \ref{basics} (\ref{Llor}).

\noindent
{\it Case A3.}
The last rule applied is ${\rm L}\limplies$.
Then the derivation ends with
$$
 \infer[{\rm L}\limplies]{
  \Gamma, S \limplies B, \Delta' \metaimplies \Sigma
  }{
  \Gamma, \Delta' \metaimplies \Sigma, S
  &
  B, \Gamma, \Delta' \metaimplies \Sigma
 } .
$$
Since
\begin{eqnarray*}
\lefteqn{(\Neg(\Gamma, \Delta') \cup \Pos(\Sigma, S)) \cap (\Nsp(\Gamma)
\cup \Pos(\Delta') \cup \Neg(\Sigma, S)) \subseteq} \\
& & 
(\Neg(\Gamma, S \limplies B, \Delta') \cup \Pos(\Sigma)) \cap (\Nsp(\Gamma)
\cup \Pos(S \limplies B, \Delta') \cup \Neg(\Sigma)) \subseteq V
\end{eqnarray*}
and
\begin{eqnarray*}
\lefteqn{(\Neg(\Gamma, B, \Delta') \cup \Pos(\Sigma)) \cap (\Nsp(\Gamma)
\cup \Pos(B, \Delta') \cup \Neg(\Sigma)) \subseteq} \\
& &
(\Neg(\Gamma, S \limplies B, \Delta') \cup \Pos(\Sigma)) \cap (\Nsp(\Gamma)
\cup \Pos(S \limplies B, \Delta') \cup \Neg(\Sigma)) \subseteq V ,
\end{eqnarray*}
we have
\[ \vdash_i \Pi_V, \Gamma, \lnot_\ast \lnot \Delta', \lnot_\ast \Sigma,
\lnot_\ast S \metaimplies \ast
\quad\mbox{and}\quad
\vdash_i \Pi_V, \Gamma, \lnot_\ast \lnot B, \lnot_\ast \lnot \Delta',
\lnot_\ast \Sigma \metaimplies \ast \]
by the induction hypothesis, and therefore, since
\[ \vdash_i \Pi_V, \Gamma, \lnot_\ast \lnot \Delta', \lnot_\ast \Sigma
\metaimplies \lnot_\ast \lnot_\ast S \]
we have
$ \vdash_i \Pi_V, \Gamma, \lnot_\ast \lnot_\ast S \limplies \lnot_\ast \lnot
B, \lnot_\ast \lnot \Delta', \lnot_\ast \Sigma \metaimplies
\lnot_\ast \lnot_\ast S , $
by ${\rm LW}$.
Thus
\[ \vdash_i \Pi_V, \Gamma, \lnot_\ast \lnot_\ast S \limplies \lnot_\ast \lnot
B, \lnot_\ast \lnot \Delta', \lnot_\ast \Sigma \metaimplies \ast
\]
by ${\rm L}\limplies$, and so
$ \vdash_i \Pi_V, \Gamma, \lnot_\ast \lnot (S \limplies B), \lnot_\ast \lnot
\Delta', \lnot_\ast \Sigma \metaimplies \ast , $
by ${\rm Cut}$ with Lemma \ref{basics} (\ref{Llimplies2}).

\noindent
{\it Case B.}
The last rule applied is an ${\rm L}$-rule, and the principal
formula is in $\Gamma$.
Since the cases for the rules ${\rm L}\land$ and ${\rm L}\lor$ are
straightforward, we review the case for the rule ${\rm L}\limplies$.

\noindent
{\it Case B1.}
The last rule applied is ${\rm L}\limplies$.
Then the derivation ends with
$$
 \infer[{\rm L}\limplies]{
  S \limplies B, \Gamma', \Delta \metaimplies \Sigma
  }{
  \Gamma', \Delta \metaimplies \Sigma, S
  &
  B, \Gamma', \Delta \metaimplies \Sigma
 } .
$$
Since
\begin{eqnarray*}
\lefteqn{(\Neg(\Gamma', \Delta) \cup \Pos(\Sigma, S)) \cap (\Nsp(\Gamma')
\cup \Pos(\Delta) \cup \Neg(\Sigma, S)) \subseteq} \\
& &
(\Neg(S \limplies B, \Gamma', \Delta) \cup \Pos(\Sigma)) \cap (\Nsp(S
\limplies B, \Gamma') \cup \Pos(\Delta) \cup \Neg(\Sigma))
\subseteq V
\end{eqnarray*}
and
\begin{eqnarray*}
\lefteqn{(\Neg(B, \Gamma', \Delta) \cup \Pos(\Sigma)) \cap (\Nsp(B, \Gamma')
\cup \Pos(\Delta) \cup \Neg(\Sigma)) \subseteq} \\
& &
(\Neg(S \limplies B, \Gamma', \Delta) \cup \Pos(\Sigma)) \cap (\Nsp(S
\limplies B, \Gamma') \cup \Pos(\Delta) \cup \Neg(\Sigma)) \subseteq V ,
\end{eqnarray*}
we have
\[ \vdash_i \Pi_V, \Gamma', \lnot_\ast \lnot \Delta, \lnot_\ast \Sigma,
\lnot_\ast S \metaimplies \ast
\quad\mbox{and}\quad
\vdash_i \Pi_V, B, \Gamma', \lnot_\ast \lnot \Delta,
\lnot_\ast \Sigma \metaimplies \ast \]
by the induction hypothesis, and therefore, since
\[ \vdash_i \Pi_V, \Gamma', \lnot_\ast \lnot \Delta, \lnot_\ast \Sigma
\metaimplies \lnot_\ast \lnot_\ast S \]
we have
$ \vdash_i \Pi_V, \lnot_\ast \lnot_\ast S \limplies \lnot_\ast \lnot_\ast B,
\Gamma', \lnot_\ast \lnot \Delta, \lnot_\ast \Sigma \metaimplies
  \lnot_\ast \lnot_\ast S , $
by ${\rm LW}$, and
\[ \vdash_i \Pi_V, \lnot_\ast \lnot_\ast B, \Gamma', \lnot_\ast \lnot \Delta,
  \lnot_\ast \Sigma \metaimplies \ast . \]
Thus
\[ \vdash_i \Pi_V, \lnot_\ast \lnot_\ast S \limplies \lnot_\ast \lnot_\ast B,
\Gamma', \lnot_\ast \lnot \Delta, \lnot_\ast \Sigma \metaimplies \ast
\]
by ${\rm L}\limplies$, and so
$ \vdash_i \Pi_V, S \limplies B, \Gamma', \lnot_\ast \lnot \Delta,
\lnot_\ast \Sigma \metaimplies \ast , $
by ${\rm Cut}$ with Lemma \ref{basics} (\ref{Llimplies1}).

\noindent
{\it Case C.}
The last rule applied is an ${\rm R}$-rule.

\noindent
{\it Case C1.}
The last rule applied is ${\rm R}\land$.
Then the derivation ends with
$$
 \infer[{\rm R}\land]{
  \Gamma, \Delta \metaimplies \Sigma', S \land S'
  }{
  \Gamma, \Delta \metaimplies \Sigma', S
  &
  \Gamma, \Delta \metaimplies \Sigma', S'
 } .
$$
Since
$ (\Neg(\Gamma, \Delta) \cup \Pos(\Sigma', S)) \cap (\Nsp(\Gamma)
\cup \Pos(\Delta) \cup \Neg(\Sigma', S)) \subseteq V $
and
$ (\Neg(\Gamma, \Delta) \cup \Pos(\Sigma', S')) \cap (\Nsp(\Gamma)
\cup \Pos(\Delta) \cup \Neg(\Sigma', S')) \subseteq V , $
we have
\[ \vdash_i \Pi_V, \Gamma, \lnot_\ast \lnot \Delta, \lnot_\ast \Sigma',
\lnot_\ast S \metaimplies \ast
\quad\mbox{and}\quad
\vdash_i \Pi_V, \Gamma, \lnot_\ast \lnot \Delta, \lnot_\ast \Sigma',
\lnot_\ast S' \metaimplies \ast \]
by the induction hypothesis, and hence
\[ \vdash_i \Pi_V, \Gamma, \lnot_\ast \lnot \Delta, \lnot_\ast \Sigma'
\metaimplies \lnot_\ast \lnot_\ast S
\quad\mbox{and}\quad
\vdash_i \Pi_V, \Gamma, \lnot_\ast \lnot \Delta, \lnot_\ast \Sigma'
\metaimplies \lnot_\ast \lnot_\ast S' . \]
Therefore
$ \vdash_i \Pi_V, \Gamma, \lnot_\ast \lnot \Delta, \lnot_\ast \Sigma'
\metaimplies \lnot_\ast \lnot_\ast (S \land S') , $
by ${\rm R}\land$ and ${\rm Cut}$ with Lemma \ref{basics}
(\ref{Rland}), and so 
$ \vdash_i \Pi_V, \Gamma, \lnot_\ast \lnot \Delta, \lnot_\ast \Sigma',
\lnot_\ast (S \land S') \metaimplies \ast . $

\noindent
{\it Case C2.}
The last rule applied is ${\rm R}\lor$.
Then the derivation ends with
$$
 \infer[{\rm R}\lor]{
  \Gamma, \Delta \metaimplies \Sigma', S \lor S'
  }{
  \Gamma, \Delta \metaimplies \Sigma', S, S'
 } .
$$
Since
$ (\Neg(\Gamma, \Delta) \cup \Pos(\Sigma', S, S')) \cap (\Nsp(\Gamma)
\cup \Pos(\Delta) \cup \Neg(\Sigma', S, S')) \subseteq V , $
we have
\[ \vdash_i \Pi_V, \Gamma, \lnot_\ast \lnot \Delta, \lnot_\ast \Sigma',
\lnot_\ast S, \lnot_\ast S' \metaimplies \ast \]
by the induction hypothesis, and hence
\[ \vdash_i \Pi_V, \Gamma, \lnot_\ast \lnot \Delta, \lnot_\ast \Sigma',
\lnot_\ast S \land \lnot_\ast S' \metaimplies \ast \]
by ${\rm L}\land$.
Therefore
$ \vdash_i \Pi_V, \Gamma, \lnot_\ast \lnot \Delta, \lnot_\ast \Sigma'
\metaimplies \lnot_\ast (\lnot_\ast S \land \lnot_\ast S') , $
and so
\[ \vdash_i \Pi_V, \Gamma, \lnot_\ast \lnot \Delta, \lnot_\ast \Sigma'
\metaimplies \lnot_\ast \lnot_\ast (S \lor S') \]
by ${\rm Cut}$ with Lemma \ref{basics} (\ref{Rlor}).
Thus
$ \vdash_i \Pi_V, \Gamma, \lnot_\ast \lnot \Delta, \lnot_\ast \Sigma',
\lnot_\ast (S \lor S') \metaimplies \ast . $

\noindent
{\it Case C3.}
The last rule applied is ${\rm R}\limplies$.
Then the derivation ends with
$$
 \infer[{\rm R}\limplies]{
  \Gamma, \Delta \metaimplies \Sigma', A \limplies S
  }{
  A, \Gamma, \Delta \metaimplies \Sigma', S
 } .
$$
Since
\begin{eqnarray*}
\lefteqn{(\Neg(\Gamma, A, \Delta) \cup \Pos(\Sigma', S)) \cap
(\Nsp(\Gamma) \cup \Pos(A, \Delta) \cup \Neg(\Sigma', S))
=} \\
& &
(\Neg(\Gamma, \Delta) \cup \Pos(\Sigma', A \limplies S)) \cap
(\Nsp(\Gamma) \cup \Pos(\Delta) \cup \Neg(\Sigma', A \limplies S))
\subseteq V ,
\end{eqnarray*}
we have
\[ \vdash_i \Pi_V, \Gamma, \lnot_\ast \lnot A, \lnot_\ast \lnot \Delta,
\lnot_\ast \Sigma', \lnot_\ast S \metaimplies \ast \]
by the induction hypothesis, and therefore, since
\[ \vdash_i \Pi_V, \Gamma, \lnot_\ast \lnot A, \lnot_\ast \lnot \Delta,
\lnot_\ast \Sigma' \metaimplies \lnot_\ast \lnot_\ast S \]
we have
$ \vdash_i \Pi_V, \Gamma, \lnot_\ast \lnot \Delta, \lnot_\ast \Sigma'
\metaimplies \lnot_\ast \lnot A \limplies \lnot_\ast \lnot_\ast S , $
by ${\rm R}\limplies$.
Thus
\[ \vdash_i \Pi_V, \Gamma, \lnot_\ast \lnot \Delta, \lnot_\ast \Sigma'
\metaimplies \lnot_\ast \lnot_\ast (A \limplies S) \]
by ${\rm Cut}$ with Lemma \ref{basics} (\ref{Rlimplies}), and so
$ \vdash_i \Pi_V, \Gamma, \lnot_\ast \lnot \Delta, \lnot_\ast
\Sigma', \lnot_\ast (A \limplies S) \metaimplies \ast . $
\qed

\begin{thm}
If
$ \vdash_c \Gamma \metaimplies A , $
then
$ \vdash_i \Pi_V, \Gamma \metaimplies A , $
where
$ V = (\Neg(\Gamma) \cup \Pos(A)) \cap (\Nsp(\Gamma) \cup \Neg(A)) . $
\end{thm}
\proof
Suppose that
$ \vdash_c \Gamma \metaimplies A , $
and let
$ V = (\Neg(\Gamma) \cup \Pos(A)) \cap (\Nsp(\Gamma) \cup \Neg(A)) . $
Then
$ \vdash_i \Pi_V, \Gamma, \lnot_\ast A \metaimplies \ast , $
by Proposition \ref{main}, and hence
\[ \vdash_i \Pi_V, \Gamma, A \limplies A \metaimplies A \]
by Lemma \ref{subst}.
Therefore
$ \vdash_i \Pi_V, \Gamma \metaimplies A . $
\qed

\begin{cor}
If
$ \vdash_c \Gamma \metaimplies A $
and
$ (\Neg(\Gamma) \cup \Pos(A)) \cap (\Nsp(\Gamma) \cup \Neg(A)) =
\emptyset , $
then
$ \vdash_i \Gamma \metaimplies A . $
\end{cor}

\section*{Acknowledgement}
The author thanks the Japan Society for the Promotion of Science
(Grant-in-Aid for Scientific Research (C) No.23540130) for partly
supporting the research.

\end{document}